# A General Family of Estimators for Estimating Population Mean in Systematic Sampling Using Auxiliary Information in the Presence of Missing Observations

[1]Manoj K. Chaudhary, [1]Sachin Malik, [2]Jayant Singh and †[1]Rajesh Singh

[1]Department of Statistics, Banaras Hindu University

Varanasi-221005, India

[2]Department of Statistics, Rajasthan University, Jaipur, India

†Corresponding author

**Abstract**

This paper proposes a general family of estimators for estimating the population mean in systematic sampling in the presence of non-response adapting the family of estimators proposed by Khoshnevisan et al. (2007). In this paper we have discussed the general properties of the proposed family including optimum property. The results have been illustrated numerically by taking an empirical population considered in the literature.

**Keywords:** Family of estimators, Auxiliary information, Mean square error, Non-response, Systematic sampling.

1. Introduction

The method of systematic sampling, first studied by Madow and Madow (1944), is used widely in surveys of finite populations. When properly applied, the methods pocks up any obvious or hidden stratification in the population and thus can be more precise than random sampling. In addition, systematic sampling is implemented easily, thus reducing costs. In this variant of random sampling, only the first unit of the sample is selected at random from the population. The subsequent units are then selected by following some definite rule.

Systematic sampling has been considered in detail by Cochran (1946) and Lahiri (1954). Reviews of the work done in the field have been given by Yates (1948) and

Buckland (1951). The application of systematic sampling to forest surveys has been illustrated by Hasel (1942), Finney (1948) and Nair and Bhargava (1951). Use of systematic sampling in estimating catch of fish has been demonstrated by Sukhatme et al. (1958). The use of auxiliary information has been permeated the important role to improve the efficiency of the estimator. Kushwaha and Singh (1989) suggested a class of almost unbiased ratio and product type estimators for estimating the population mean using jack-knife technique initiated by Quenouille (1956). Afterward Banarasi et al. (1993) and Singh and Singh (1998) have proposed the estimators of population mean using auxiliary information in systematic sampling.

Khoshnevisan et al. (2007) suggested a general family of estimators for estimating the populations mean using known values of some population parameters in simple random sampling, given by

$$t = \bar{y}\left[\frac{a\bar{X}+b}{\alpha(a\bar{x}+b)+(1-\alpha)(a\bar{X}+b)}\right]^g \qquad (1.1)$$

where $\bar{y}$ and $\bar{x}$ are the sample means of study and auxiliary variables respectively. $\bar{X}$ is the population mean of auxiliary variable $X$. $a \neq 0$ and $b$ are either real numbers or functions of known parameters of auxiliary variable. $\alpha$ and $g$ are the real numbers which are to be determined. Here we would like to mention that the choice of the estimator depends on the availability and values of the various parameter(s) used (for choice of the parameters a and b refer to Singh et al. (2008) and Singh and Kumar(2011)).

In this paper we have proposed a general family of estimators for estimating the population mean in systematic sampling using auxiliary information in the presence of non-response following Khoshnevisan et al. (2007). We have also derived the expressions for minimum mean square errors (MSE) of the family with respect to $\alpha$. A comparative study is also carried out to compare the optimum estimators of the family with respect to usual mean estimator with the help of numerical data.

## 2. Proposed Family of Estimators

Let us suppose that a population consists of $N$ units numbered from 1 to $N$ in some order and a sample of size $n$ is to be drawn such that $N = nk$ ($k$ is an integer). Thus there will be k samples each of n units and we select one sample from the set of k samples. Let Y and X be the study and auxiliary variable with respective means $\bar{Y}$ and $\bar{X}$. Let us consider $y_{ij}(x_{ij})$ be the $j^{th}$ observation in the $i^{th}$ systematic sample under study (auxiliary) variable ($i=1...k : j=1...n$).

Wwe assume that the non-response is observed only on study variable and auxiliary variable is free from non-response. Using Hansen-Hurwitz (1946) technique of sub-sampling of non-respondents, the estimator of population mean $\bar{Y}$, can be defined as

$$\bar{y}^* = \frac{n_1 \bar{y}_{n1} + n_2 \bar{y}_{h_2}}{n} \qquad (2.1)$$

where $\bar{y}_{n1}$ and $\bar{y}_{h_2}$ are, respectively the means based on $n_1$ respondent units from the systematic sample of n units and sub-sample of $h_2$ units selected from $n_2$ non-respondent units in the systematic sample. The estimator of population mean $\bar{X}$ of auxiliary variable based on the systematic sample of size n units, is given by

$$\bar{x} = \frac{1}{n} \sum_{j=1}^{n} x_{ij} \qquad (i=1...k) \qquad (2.2)$$

Obviously, $\bar{y}^*$ and $\bar{x}$ are unbiased estimators. The variance expression for the estimators $\bar{y}^*$ and $\bar{x}$ are, respectively, given by

$$V(\bar{y}^*) = \frac{N-1}{nN}\{1+(n-1)\rho_Y\}S_Y^2 + \frac{L-1}{n}W_2 S_{Y2}^2 \qquad (2.3)$$

and $\qquad V(\bar{x}) = \frac{N-1}{nN}\{1+(n-1)\rho_X\}S_X^2 \qquad (2.4)$

where $\rho_Y$ and $\rho_X$ are the correlation coefficients between a pair of units within the systematic sample for the study and auxiliary variables respectively. $S_Y^2$ and $S_X^2$ are respectively the mean squares of the entire group for study and auxiliary variable. $S_{Y2}^2$ be the mean square of non-response group under study variable, $W_2$ is the non-response rate in the population and $L = \frac{n_2}{h_2}$.

Let us assume that the population mean $\overline{X}$ is known. Thus the usual ratio and product estimators of the population mean $\overline{Y}$ under non-response based on a systematic sample of size $n$, can be respectively defined as

$$\overline{y}_R^* = \frac{\overline{y}^*}{\overline{x}}\overline{X} \qquad (2.5)$$

and

$$\overline{y}_P^* = \frac{\overline{y}^*\overline{x}}{\overline{X}} \qquad (2.6)$$

To obtain the biases and mean square errors, we use large sample approximation.

$$\overline{y}^* = \overline{Y}(1+e_0)$$

$$\overline{x} = \overline{X}(1+e_1)$$

such that $E(e_0) = E(e_1) = 0$, and

$$E(e_0^2) = \frac{V(\overline{y}^*)}{\overline{Y}^2} = \frac{N-1}{nN}\{1+(n-1)\rho_Y\}C_Y^2 + \frac{L-1}{n}W_2\frac{S_{Y2}^2}{\overline{Y}^2},$$

$$E(e_1^2) = \frac{V(\overline{x})}{\overline{X}^2} = \frac{N-1}{nN}\{1+(n-1)\rho_X\}C_X^2,$$

and $E(e_0 e_1) = \dfrac{Cov(\overline{y}^*,\overline{x})}{\overline{Y}\overline{X}} = \dfrac{N-1}{nN}\{1+(n-1)\rho_Y\}^{1/2}\{1+(n-1)\rho_X\}^{1/2}\rho C_Y C_X$

where $C_Y$ and $C_X$ are the coefficients of variation of study and auxiliary variables respectively.

Expressing the equations (2.5) and (2.6) in terms of $e_i$'s $(i=0,1)$ and taking expectations the bias expressions of the estimators of $\overline{y}_R^*$ and $\overline{y}_P^*$, are respectively given by

$$B(\overline{y}_R^*) = \frac{N-1}{nN}\overline{Y}\{1+(n-1)\rho_X\}(1-K\rho^*)C_X^2 \qquad (2.7)$$

and

$$B(\overline{y}_P^*) = \frac{N-1}{nN}\overline{Y}\{1+(n-1)\rho_X\}K\rho^* C_X^2 \qquad (2.8)$$

where, $\rho^* = \dfrac{\{1+(n-1)\rho_Y\}^{1/2}}{\{1+(n-1)\rho_X\}^{1/2}}$ and $K = \rho \dfrac{C_Y}{C_X}$.

The mean square errors (MSE's) of $\overline{y}_R^*$ and $\overline{y}_P^*$, are respectively, given by

$$\mathrm{MSE}\left(\overline{y}_R^*\right) = \dfrac{N-1}{nN}\overline{Y}^2\{1+(n-1)\rho_X\}\left[\rho^{*2}C_Y^2 + \left(1-2K\rho^*\right)C_X^2\right] + \dfrac{L-1}{n}W_2 S_{Y2}^2$$

(2.9)

and

$$\mathrm{MSE}\left(\overline{y}_P^*\right) = \dfrac{N-1}{nN}\overline{Y}^2\{1+(n-1)\rho_X\}\left[\rho^{*2}C_Y^2 + \left(1+2K\rho^*\right)C_X^2\right] + \dfrac{L-1}{n}W_2 S_{Y2}^2$$

(2.10)

Motivated by Khoshnevisan et al. (2007), we now define a family of estimators of population mean $\overline{Y}$ based on a systematic sample of size $n$ in the presence of non-response as

$$t^* = \overline{y}^*\left[\dfrac{a\overline{X}+b}{\alpha(a\overline{x}+b)+(1-\alpha)(a\overline{X}+b)}\right]^g$$

(2.11)

This family can generate the non-response versions of a number of estimators of population mean $\overline{Y}$ including the usual ratio and product estimators on different choices of $a$, $b$, $\alpha$ and $g$.

### 2.1 Properties of $t^*$

Expressing $t^*$ in terms of $e_i$'s, we get

$$t^* = \overline{y}(1+e_0)(1+\alpha\lambda e_1)^{-g}$$

(2.12)

where $\lambda = \dfrac{a\overline{X}}{a\overline{X}+b}$.

We assume that $|\lambda e_1| < 1$ so that the right–hand side of the equation (2.12) is expandable in terms of power series. Expanding the right–hand side of the equation (2.12) and neglecting the terms in $e_i$'s having power greater than two, we have

$$t^* - \overline{Y} = \overline{Y}\left[e_0 - g\alpha\lambda e_1 + \frac{g(g+1)}{2}\alpha^2\lambda^2 e_1^2 - g\alpha\lambda e_0 e_1\right] \quad (2.13)$$

Taking expectation both the sides of equation (2.13), we get the bias of $t^*$ up to the first order of approximation, as

$$B(t^*) = \frac{N-1}{nN}\overline{Y}\{1+(n-1)\rho_X\}C_X^2\left[\frac{g(g+1)}{2}\alpha^2\lambda^2 - g\alpha\lambda K\rho^*\right] \quad (2.14)$$

Squaring both the sides of the equation (2.13) and then taking the expectation, we obtain the MSE of $t^*$ up to the first order of approximation, as

$$MSE(t^*) = \frac{N-1}{nN}\overline{Y}^2\{1+(n-1)\rho_X\}\left[\rho^{*2}C_Y^2 + \left(g^2\alpha^2\lambda^2 - 2g\alpha\lambda\rho^* K\right)C_X^2\right] + \frac{(L-1)}{n}W_2 S_{Y2}^2$$

(2.15)

**2.2 Optimum Choice of $\alpha$**

In order to obtain the minimum MSE of $t^*$, we differentiate the MSE of $t^*$ with respect to $\alpha$ and equating the derivative to zero, we get

$$\frac{N-1}{nN}\overline{Y}^2\{1+(n-1)\rho_X\}\left[\alpha g^2\lambda^2 - g\lambda\rho^* K\right]C_X^2 = 0 \quad (2.16)$$

The equation (2.16) provides the optimum values of $\alpha$ as

$$\alpha = \frac{\rho^* K}{g\lambda} \quad (2.17)$$

Putting the optimum value of $\alpha$ from equation (2.17) into the equation (2.15), we get the minimum MSE of $t^*$, as

$$MSE(t^*)_{min} = \frac{N-1}{nN}\overline{Y}^2\{1+(n-1)\rho_X\}\left[C_Y^2 - K^2 C_X^2\right]\rho^{*2} + \frac{(L-1)}{n}W_2 S_{Y2}^2 \quad (2.18)$$

The minimum MSE of $t^*$, is same as the mean square error of the usual regression estimator in systematic sampling under non-response.

## 3. Empirical Study

For numerical illustration, we have considered the data given in Murthy (1967, p. 131-132). The data are based on length (X) and timber volume (Y) for 176 forest strips. Murthy (1967, p.149) and Kushwaha and Singh (1989) reported the values of intraclass correlation coefficients $\rho_X$ and $\rho_Y$ approximately equal for the systematic sample of size 16 by enumerating all possible systematic samples after arranging the data in ascending order of strip length. The details of population parameters are :

$N = 176$, $\quad n = 16$, $\quad \overline{Y} = 282.6136$, $\quad \overline{X} = 6.9943$,

$S_Y^2 = 24114.6700$, $\quad S_X^2 = 8.7600$, $\quad \rho = 0.8710$,

$S_{Y2}^2 = \dfrac{3}{4} S_Y^2 = 18086.0025$.

Table 1 shows the percentage relative efficiency (PRE) of $t^*$ (optimum) with respect to $\overline{y}^*$ for the different choices of $W_2$ and $L$.

**Table 1**: PRE of $t^*$ (optimum) with respect to $\overline{y}^*$

| $W_2$ | L | PRE |
|---|---|---|
| 0.1 | 2.5 | 407.48 |
|  | 2.5 | 404.18 |
|  | 3.0 | 400.94 |
|  | 3.5 | 397.77 |
| 0.2 | 2.5 | 400.94 |
|  | 2.5 | 394.67 |
|  | 3.0 | 388.66 |
|     | 3.5 | 382.89 |
| 0.3 | 2.5 | 394.67 |
|     | 2.5 | 385.74 |
|     | 3.0 | 377.34 |
|     | 3.5 | 369.42 |
| 0.4 | 2.5 | 403.22 |
|     | 2.5 | 377.34 |
|     | 3.0 | 366.88 |
|     | 3.5 | 357.17 |

## 4. Conclusion

In this paper, we have proposed a general family of estimators of population mean in systematic sampling using an auxiliary variable in the presence of non-response. The optimum property of the family has been discussed. The study concludes that the suggested family converges to the usual regression estimator of population mean in systematic sampling under non-response if the parameter $\alpha$ attains its optimum value. From Table 1, it can easily be seen that the estimator $t^*$ (optimum) performs always better than the usual estimator $\bar{y}^*$. It is also observed that the percentage relative efficiency (PRE) of $t^*$ (optimum) with respect to $\bar{y}^*$ decreases with increase in non-response rate $W_2$ as well as $L$.